# HARMONIC CONTINUOUS-TIME BRANCHING MOMENTS

By Didier Piau

Université Lyon 1 and Université Grenoble 1

We show that the mean inverse populations of nondecreasing, square integrable, continuous-time branching processes decrease to zero like the inverse of their mean population if and only if the initial population $k$ is greater than a first threshold $m_1 \geq 1$. If, furthermore, $k$ is greater than a second threshold $m_2 \geq m_1$, the normalized mean inverse population is at most $1/(k - m_2)$. We express $m_1$ and $m_2$ as explicit functionals of the reproducing distribution, we discuss some analogues for discrete time branching processes and link these results to the behavior of random products involving i.i.d. nonnegative sums.

**Introduction.** We consider nondecreasing continuous-time branching processes $(z_t)_{t \geq 0}$ with initial population $z_0 := k \geq 1$, split intensity $s(t)$ and offspring distribution $p(t) * \delta_1$ at time $t \geq 0$, where

$$p(t) = \sum_{i \geq 0} p(i,t)\, \delta_i.$$

In other words, $z_t$ is the size at time $t$ of a population which evolves as follows. During the interval $[t, t+dt)$, each individual living at time $t$ survives with probability $1 - s(t)\, dt + o(dt)$. For every $i \geq 0$, the individual dies and is instantaneously replaced by $i + 1$ individuals, independently of the behavior of the other individuals living at time $t$, with probability $p(i,t)s(t)\, dt + o(dt)$. The process $(z_t)_t$ is pathwise nondecreasing, hence $z_t \geq k$ with full probability and the harmonic moments $\mathbb{E}_k\{1/z_t\}$ are well defined.

Harmonic moments of branching processes are needed to build unbiased estimators of the offspring distributions from samples of Markov branching processes in noncanonical situations, for instance when the state of an individual depends on the number of its siblings. Examples of this situation in a discrete-time setting arise in the context of applications to molecular biology; see [8] for a presentation and [5, 6, 7] for a mathematical study. Also, harmonic moments reflect mainly lower deviations of the branching process and may be viewed as an integrated way to quantify these deviations; see [3, 4] and the thorough exposition in [2].

In this paper, we provide sharp bounds for $\mathbb{E}_k\{1/z_t\}$ if and when $\mathbb{E}_k\{1/z_t\}$ decreases to zero roughly like $1/\mathbb{E}_k\{z_t\}$ when $t \to \infty$. The emphasis is on explicit and computable upper bounds and the setting is restricted to square integrable offspring distributions.

**1. Main results.** This section is composed as follows. We first introduce some definitions (Section 1.1) and explain a canonical reduction of the model (Section 1.2). Then we explain our main results the continuous-time case (Section 1.3) and the discrete-time case (Section 1.4). Finally, we sketch some relations of our results to previous work (Section 1.5) and describe the overall plan of the rest of the paper (Section 1.6).

1.1. *Notation.*

DEFINITION 1.1.  For every $t \geq 0$, $L(t) \geq 0$ denotes a random variable of distribution $p(t)$ and

$$\mathbb{E}\{L(t)\} = \sum_{i \geq 0} i p(i, t), \qquad M(t) := \int_0^t \mathbb{E}\{L(u)\} s(u)\, du.$$

The function $h_k$ and the quantity $h(k)$ are defined by

$$h_k(t) := e^{M(t)}\, \mathbb{E}_k\{1/z_t\}, \qquad h(k) := \sup_{t \geq 0} h_k(t).$$

From elementary computations,

$$\mathbb{E}_k\{z_t\} = k e^{M(t)}.$$

Since the function $z \mapsto 1/z$ is convex, $h_k(\cdot)$ is nondecreasing, $h(k)$ is the limit of $h_k(\cdot)$ at infinity, $h(k) \geq 1/k = h_k(0)$ and, for every $t \geq 0$,

$$e^{-M(t)}/k \leq \mathbb{E}_k\{1/z_t\} \leq h(k) e^{-M(t)}.$$

Hereafter, our aim is to provide sharp upper and lower bounds of $h(k)$.

1.2. *Reduction.*  In this section, we show that one can assume $s(t) = 1$ and $L(t) \geq 1$ for every $t \geq 0$.

DEFINITION 1.2.  For every random variable $L$ such that $\mathbb{P}(L = 0) \neq 1$, $L'$ denotes a random variable distributed like $L$ and conditioned by $\{L \neq 0\}$.

DEFINITION 1.3. Let $(*)$ denote the condition that the function $\sigma(\cdot)$ is unbounded, where, for every $t \geq 0$,

$$\sigma(t) := \int_0^t (1 - p(0,u))s(u)\, du.$$

The individuals who survive and those who die and are replaced by exactly one individual have the same net effect on the overall count $z_t$ of the population. This proves Proposition 1.4.

PROPOSITION 1.4 (Canonical reduction). *Let $(z'_t)_t$ denote the branching process of parameters $p'(\cdot)$ and $s'(\cdot) \equiv 1$, where $p'(t)$ denotes the distribution of $L'(\sigma^{-1}(t))$. Then the distributions of $(z_t)_t$ and $(z'_{\sigma(t)})_t$ coincide. As a consequence, provided the condition $(*)$ is met, the quantities $h(k)$ for $(z_t)_t$ and for $(z'_t)_t$ coincide.*

Hence, under condition $(*)$, one can assume without loss of generality (and we do assume this from now on) that $L(t) \geq 1$ for every $t$, almost surely, and that $s(\cdot) \equiv 1$.

1.3. *Results in continuous time.* The case when the random tree associated to the branching process is regular is, in a sense, extremal. We first deal with this case.

THEOREM 1 (Degenerate case). *Assume that $p(t) = \delta_i$ for every $t \geq 0$, with $i \geq 1$. Then, for every $k \leq i$, $h(k)$ is infinite and for every $k \geq i+1$,*

$$h(k) = 1/(k-i).$$

In the general case, the parameters $\mathfrak{m}_1$ and $\mathfrak{m}_2$ defined below play a crucial role.

DEFINITION 1.5. For every integer-valued, not identically zero, square integrable $L \geq 0$, let

$$\mathfrak{m}_1\{L\} := \mathbb{E}\{L \,|\, L \neq 0\}, \qquad \mathfrak{m}_2\{L\} := \mathbb{E}\{L^2\}/\mathbb{E}\{L\}.$$

For every $L$, $\mathfrak{m}_1\{L\} = \mathfrak{m}_1\{L'\}$ and $\mathfrak{m}_2\{L\} = \mathfrak{m}_2\{L'\}$ (see Definition 1.2), and

$$1 \leq \mathfrak{m}_1\{L\} \leq \mathfrak{m}_2\{L\}.$$

If $L$ assumes more than one nonzero value, these inequalities are strict. Our next result deals with the case when $p(t)$ does not depend on $t$:

THEOREM 2 (Homogeneous case). *Assume that $L(t)$ is distributed as $L \geq 0$, for every positive $t$, and that $L$ is square integrable and not identically zero.*

(a) *For every $k \leq \mathfrak{m}_1\{L\}$, $h(k)$ is infinite. In particular, $h(1)$ is infinite.*
(b) *For every $k > \mathfrak{m}_1\{L\}$, $h(k)$ is finite and*
$$h(k) > 1/(k - \mathfrak{m}_1\{L\}).$$
(c) *For every $k > \mathfrak{m}_2\{L\}$,*
$$1/(k - \mathfrak{m}_1\{L\}) < h(k) \leq 1/(k - \mathfrak{m}_2\{L\}).$$

Finally, we state our result in the general case:

THEOREM 3 (General case). *Assume that there exist finite positive constants $\mathfrak{m}_1^\pm$ and $\mathfrak{m}_2^+$ such that, for every $t$,*
$$\mathfrak{m}_1^- \leq \mathfrak{m}_1\{L(t)\} \leq \mathfrak{m}_1^+, \qquad \mathfrak{m}_2\{L(t)\} \leq \mathfrak{m}_2^+.$$
*Then the following hold:*

(a) *For every $k \leq \mathfrak{m}_1^-$, $h(k)$ is infinite.*
(b) *For every $k > \mathfrak{m}_1^+$, $h(k)$ is finite.*
(c) *For every $k > \mathfrak{m}_2^+$, $1/(k - \mathfrak{m}_1^-) \leq h(k) \leq 1/(k - \mathfrak{m}_2^+)$.*

We explain below how to deduce from Theorems 2 and 3 upper bounds of $h(k)$ in the ranges $\mathfrak{m}_1\{L\} < k \leq \mathfrak{m}_2\{L\}$ and $\mathfrak{m}_1^+ < k \leq \mathfrak{m}_2^+$, respectively. Finally, note that in Theorem 3, $\mathfrak{m}_1\{L(t)\}$ and $\mathfrak{m}_2\{L(t)\}$ can be independent on $t$, even when the distribution of $L(t)$ depends on $t$.

1.4. *Results in discrete time.* This section refines results in [7] which we recall in Section 2.1, and uses notation from Definition 2.4 in Section 2.1.

THEOREM 4. *In the setting of Theorem 2, for every $k > \mathfrak{m}_2\{L\}$,*
$$H(k) \leq 1/(k - \mathfrak{m}_2\{L\}).$$
*In the setting of Theorem 3, for every $k > \mathfrak{m}_2^+$,*
$$H(k) \leq 1/(k - \mathfrak{m}_2^+).$$

Overall, the situation is less clear in discrete time than in continuous time; note in particular that Theorem 4 does not yield the exact set of populations $k$ such that $H(k)$ is finite.

1.5. *Relation to previous work.* In the discrete-time setting, this paper deals with offspring distributions $\delta_1 * p$ such that $p(0) \neq 0$, which belong to the so-called Schröder's case; see [2], for instance. With the notation of our paper, Schröder's exponent $\alpha$ is defined by the relation

$$p(0)(1 + \mathbb{E}\{L\})^\alpha = 1.$$

When the distribution $p$ is $L \log L$ integrable and when $Z(0) = 1$, the martingale $Z(n)/(1 + \mathbb{E}\{L\})^n$ converges almost surely to a random variable $W$, which is almost surely positive and finite. Introducing i.i.d. copies $(W_i)_{i \geq 1}$ of $W$ and using the fact that $(1 + \mathbb{E}\{L\})^n/Z(n)$ is a positive submartingale, one gets

$$H(k) = \mathbb{E}(1/(W_1 + \cdots + W_k)).$$

Likewise, in our continuous-time setting, when $z_0 = 1$, $z_t/e^{M(t)}$ converges almost surely to $w$, say, and using i.i.d. copies $(w_i)_{i \geq 1}$ of $w$, one gets

$$h(k) = \mathbb{E}(1/(w_1 + \cdots + w_k)).$$

Returning to the discrete-time setting, $W$ has distribution $\omega(x)\,dx$, the function $\omega(\cdot)$ is continuous (see [1]), and there exist positive and finite constants $c_1$ and $c_2$ such that, when $x \to 0$,

$$c_1 x^{1-\alpha} \leq \omega(x) \leq c_2 x^{1-\alpha}.$$

Using this, one can check that $\mathbb{E}_1(1/W^\beta)$ is finite if and only if

$$\alpha > \beta.$$

Likewise, starting from $Z_0 = k$, $W$ is distributed like the sum of $k$ i.i.d. copies of $W$ starting from one individual. Hence, $H(k)$ is finite if and only if

$$k\alpha > 1.$$

One sees that, in a loose sense, the initial population $k$ plays the rôle of $1/\beta$. For instance, by convexity, it is a quite general and easily-demonstrated fact that if $\mathbb{E}_1(1/W^{1/k})$ is finite, then $H(k) = \mathbb{E}_k(1/W)$ is also finite. For more remarks on this, see Section 1.5 of [7].

1.6. *Plan.* The rest of the paper is organized as follows. In Section 2, we recall some results from [7] about the harmonic moments of branching processes in discrete time. We then refine these results and show how to apply them to the continuous-time setting. This section culminates with the statement of Theorem 5, which we prove later. In Section 3, we prove Theorem 1 and state some additional remarks. Section 4 presents a stochastic recursion, basic for our study. This allows the exact determination of the initial populations $k$ such that $h(k)$ is finite. In Section 5, we prove part

(b) of Theorem 2. In Section 6, we prove Theorem 5, stated in Section 2. In Section 7, we use the upper bounds for initial populations above $\mathfrak{m}_2$ to deduce upper bounds for smaller initial populations, namely between $\mathfrak{m}_1$ and $\mathfrak{m}_2$. Finally, we explain in Section 8 how our results on branching processes yield the limits of specific random products, which involve some sums of nonnegative i.i.d. random variables.

## 2. Continuous-time and discrete-time harmonic moments.

2.1. *Discrete-time harmonic moments.* We consider branching processes in discrete time $(Z(n))_{n \geq 0}$ with reproducing distribution $p * \delta_1$. Hence, introducing i.i.d. random variables $L(n,i) \geq 1$ with distribution $p * \delta_1$, we set

$$Z(n+1) := \sum_{i=1}^{Z(n)} L(n,i).$$

As is well known, continuous-time branching processes are limits of discrete-time ones. To see this, for every distribution $p$, we introduce a family $(Z_u)_u$ of branching processes in discrete time, based on $p$ and defined as follows:

DEFINITION 2.1. For each $u$ in $(0,1]$, $L_u \geq 0$ denotes a random variable of distribution $p_u$ given by

(1) $$p_u := (1-u)\delta_0 + up.$$

Let $(Z_u(n))_{n \geq 0}$ denote a nondecreasing, integrable Galton–Watson process, indexed by the nonnegative integers, starting from $Z_u(0) := k$ and with reproducing distribution $\delta_1 * p_u$. Then, when $u \to 0$, $Z_u$ converges to $(z_t)_t$ in distribution, at least in the sense that, for each fixed $t$, $Z_u(\lfloor t/u \rfloor)$ converges to $z_t$ in distribution. Thus, results about the harmonic moments of the discrete-time branching processes $Z_u$, when suitably uniform over the parameter $u$, yield information about the harmonic moments of the continuous-time branching process $(z_t)_t$. To develop this strategy, we borrow some definitions and results from [7], with some modifications.

DEFINITION 2.2 ([7]). For every nonnegative $m$, $g_{0,m}(s) := s^m$. For every nonnegative $m$ and every positive $c$,

$$g_{c,m}(s) := (1 + m(1-s^c))^{-1/c}.$$

For every nonnegative $c$, $\mathcal{G}_c$ denotes the space of nonnegative random variables $L$ such that, for every $s$ in $[0,1]$, $\mathbb{E}\{s^L\} \leq g_{c,m}(s)$ with $m := \mathbb{E}\{L\}$.

We recall the following facts and remarks from [7]. First,
$$g_{c,m}(s) =: \mathbb{E}\{s^{L_{c,m}}\}$$
is the generating function of a nonnegative random variable $L_{c,m}$, every moment of which is finite, such that $\mathbb{E}\{L_{c,m}\} = m$, and whose distribution is either $\delta_m$ if $c = 0$, or is supported by the set of nonnegative multiples of $c$ if $c$ is positive (viz., the distribution of $L_{c,m}/c$ is negative binomial). We stress that one compares the generating functions of $L$ and of $L_{c,m}$, two random variables whose means are equal, but that, the case when $c$ is an integer being excluded, $L_{c,m}$ is not an integer-valued random variable. Hence there does not exist, in general, any branching process based on the distribution of $L_{c,m}$.

Second, every $L$ in $\mathcal{G}_c$ is square integrable and its variance $\sigma^2\{L\}$ is at most $\sigma^2\{L_{c,m}\} = cm(m+1)$. Finally, the family $(\mathcal{G}_c)_{c \geq 0}$ is increasing and every square integrable $L$ belongs to $\mathcal{G}_c$ if $c$ is sufficiently large. Using Definition 2.3 below, we rephrase this as Proposition 2.6.

DEFINITION 2.3. For every integer-valued random variable $L \geq 0$, let
$$\mathfrak{C}\{L\} := \inf\{c \geq 0\, ; L \in \mathcal{G}_c\}.$$
When the context is clear, we use $\mathfrak{C}$ for $\mathfrak{C}\{L\}$.

Our next definition is the analogue in discrete time of Definition 1.1:

DEFINITION 2.4. For every integer-valued random variable $L \geq 0$, let
$$H_n(k) := (1 + \mathbb{E}\{L\})^n \mathbb{E}_k\{1/Z(n)\}, \qquad H(k) := \lim_{n \to \infty} H_n(k) = \sup_n H_n(k).$$

PROPOSITION 2.5 ([7]). *For every integer $k > \mathfrak{C}\{L\}$,*
$$H(k) \leq 1/(k - \mathfrak{C}\{L\}).$$

Proposition 2.6 below shows that the result of Proposition 2.5 is not empty since $\mathfrak{C}\{L\}$ is finite, at least for (in fact, exactly for) square integrable random variables $L$.

PROPOSITION 2.6 ([7]). *For every non-square-integrable $L$, $\mathfrak{C}\{L\}$ is infinite. For every square integrable $L$, $\mathfrak{C}\{L\}$ is finite and*
$$\mathfrak{C}\{L\} \geq \frac{\sigma^2\{L\}}{\mathbb{E}\{L\}(\mathbb{E}\{L\} + 1)}.$$

2.2. *From discrete time to continuous time.* We apply the results of the preceding section to continuous-time branching processes. We recall that Definition 2.1 in Section 2.1 introduces random variables $L_u$ such that $\mathbb{E}\{L_u\} = 1 + u\,\mathbb{E}\{L\}$.

DEFINITION 2.7. Let
$$H_n^u(k) := (1+u\mathbb{E}\{L\})^n\, \mathbb{E}_k\{1/Z_u(n)\}, \qquad H^u(k) := \lim_{n\to\infty} H_n^u(k) = \sup_{n\geq 0} H_n^u(k).$$

As a direct application of results in discrete time to our setting, assume that $p_u$ belongs to $\mathcal{G}_{c(u)}$ for every sufficiently small $u$. Then, for every $n$, every sufficiently small $u$ and every $k > c(u)$, Proposition 2.5 in Section 2.1 yields
$$H_n^u(k) \leq (1 + u\mathbb{E}\{L\})^n \mathbb{E}_k\{1/(Z_u(n) - c(u))\} \leq 1/(k - c(u)).$$
Hence, $H_{\lfloor t/u \rfloor}^u(k) \leq 1/(k - c(u))$ and
$$h_k(t) \leq \liminf_u 1/(k - c(u)).$$
This motivates Definition 2.8 and proves Proposition 2.9 below.

DEFINITION 2.8. For every integer-valued random variable $L \geq 0$, let
$$\mathfrak{C}_0\{L\} := \sup_{0 < u \leq 1} \mathfrak{C}\{L_u\}.$$
When the context is clear, we use $\mathfrak{C}_0$ for $\mathfrak{C}_0\{L\}$.

PROPOSITION 2.9. *For every integer $k > \mathfrak{C}_0\{L\}$,*
$$h(k) \leq 1/(k - \mathfrak{C}_0\{L\}).$$
*Hence, for every nonnegative $t$,*
$$\mathbb{E}_k\{1/z_t\} \leq \mathbb{E}_k\{1/(z_t - \mathfrak{C}_0\{L\})\} \leq e^{-mt}/(k - \mathfrak{C}_0\{L\}).$$

2.3. *Identification of $\mathfrak{C}_0$.* The main drawback of the results in Sections 2.1 and 2.2 is that the parameters $\mathfrak{C}\{L_u\}$ (and hence, presumably, the value of $\mathfrak{C}_0\{L\}$) are virtually unknown. In fact, while the exact value of $\mathfrak{C}\{L\}$, for a given random variable $L$, may indeed be difficult to obtain, it transpires that the value of $\mathfrak{C}_0\{L\}$, which is the only one relevant in the continuous-time setting, is quite simple.

Approaching this result, we first note that $\mathfrak{C}\{L_u\}$ is a monotone function of the parameter $u$ in $(0, 1]$.

PROPOSITION 2.10. *For every c, u and L, if $L_u$ belongs to $\mathcal{G}_c$, then L belongs to $\mathcal{G}_c$. In other words, the function $u \mapsto \mathfrak{C}\{L_u\}$ is nonincreasing for u in $(0, 1]$. As a consequence, the parameter $\mathfrak{C}_0\{L\}$, defined as a supremum, is also given by*

$$\mathfrak{C}_0\{L\} = \lim_{u \to 0} \mathfrak{C}\{L_u\}.$$

Hence $\mathfrak{C}\{L\}$ describes p and $\mathfrak{C}_0\{L\}$ describes the limit of $p_u$ when $u \to 0$.

Our next result is anecdotal and stated for the sake of completeness. It includes a generalization of Definition 2.4.

PROPOSITION 2.11. *For any positive integer k and any real number $b < k$, let*

$$h_k(t, b) := e^{mt} \mathbb{E}_k\{1/(z_t - b)\}, \qquad h(k, b) := \sup_{t \geq 0} h_k(t, b).$$

*The sequences $(h(k))_{k \geq 1}$ and $(h(k, b))_{k > b}$ are then convex.*

Our main result in this section provides some tight bounds for $\mathfrak{C}\{L\}$ and the exact value of $\mathfrak{C}_0\{L\}$, for every square integrable L.

THEOREM 5. *For every integer-valued, nonnegative, square integrable L,*

$$\mathfrak{m}_2\{L\} - \frac{\mathbb{E}\{L\}}{1 + \mathbb{E}\{L\}} \leq \mathfrak{C}\{L\} \leq \mathfrak{m}_2\{L\}.$$

*As a consequence,*

$$\mathfrak{C}_0\{L\} = \mathfrak{m}_2\{L\}.$$

The upper and the lower bounds that the first part of Theorem 5 yields for $\mathfrak{C}\{L_u\}$ both converge to $\mathfrak{m}_2\{L\}$ when $u \to 0$, hence the second part is a direct consequence of the first part. Finally, Theorem 5 and Proposition 2.9 imply part (c) of Theorems 2 and 3. Together with Proposition 2.10, they also yield Theorem 4. We prove Theorem 5 in Section 6.

2.4. *Proofs.*

PROOF OF PROPOSITION 2.10. Since $\mathbb{E}\{s^{L_u}\} = 1 - u + u\,\mathbb{E}\{s^L\}$, it is enough to show that $g_{c,um} \leq 1 - u + u\,g_{c,m}$. In turn, this follows from the convexity of the function $m \mapsto g_{c,m}(s)$ for every fixed s. $\square$

PROOF OF PROPOSITION 2.11. This is elementary. For every positive i.i.d. $\zeta$ and $\zeta_i$ and every $k \geq 1$, let $\sigma_k := \zeta_1 + \cdots + \zeta_k$ and $\eta_k := \mathbb{E}\{1/\sigma_k\}$. Then the sequence $(\eta_k)_k$ is convex. To see this, note that the convexity is equivalent to $(\eta_k - \eta_{k+1})_k$ being nonincreasing and that this holds true because

$$\eta_k - \eta_{k+1} = \mathbb{E}\{\zeta/((\zeta + \sigma_k)\sigma_k)\}$$

and because $(\sigma_k)_k$ is nondecreasing. Applying this convexity to $Z_u(n)$ and considering the limit when $n \to \infty$ shows that $H^u$ is a convex sequence for every $u$. Taking the limit when $u \to 0$ shows that the sequence $h$ is convex.

The proof for $h(k, b)$ is similar and hence omitted. □

**3. Regular trees.** The simplest case is when the trees are regular, that is, when the branching process is ruled by a distribution $p := \delta_i$ with $i \geq 1$. We already mentioned that the case $(1 - \lambda)\delta_0 + \lambda\delta_i$ with $\lambda$ in $(0, 1]$ is equivalent, since this corresponds to a change of time of the branching mechanism. The associated random tree is regular with degree $i + 1$ and the length of every edge is exponential of mean $1/\lambda$. The distribution $p_u$ is Bernoulli 0 or $i$ with $p_u(i) = u\lambda$. Such Bernoulli distributions belong to $\mathcal{G}_i$. In fact (see [7]),

$$\mathfrak{C}_0 = \mathfrak{m}_1 = \mathfrak{m}_2 = i.$$

Hence, for every $k \geq i + 1$,

$$h(k) \leq 1/(k - i)$$

and this method cannot yield a finite upper bound in the $k = i$ case. This is for good reason since, as we show below, $h(i)$ is in fact infinite.

The distribution of $z_t$ when $z_0 = 1$ is an affine transform of the negative binomial distribution of index $i$ and mean $e^{\lambda t}$. More precisely, for every $|v| \leq 1$,

$$\mathbb{E}_1\{v^{z_t}\} = \frac{v}{(e^{\lambda t} - (e^{\lambda t} - 1)v^i)^{1/i}}.$$

Hence, for every positive $t$,

$$h_i(t) = \frac{\lambda t e^{\lambda t}}{i(e^{\lambda t} - 1)}$$

and $h_i$ increases from $h_i(0) = 1/i$ to $h(i) = \infty$. This shows that $h_k$ is unbounded when $k \leq i$. On the other hand, for every $k \geq i + 1$, the result above is sharp since (see below) the function $h_k$ increases from $h_k(0) = 1/k$ to

$$h(k) = 1/(k - i).$$

For branching processes in discrete time (see [7]), the true value of $\mathfrak{C}\{L\}$ is

$$\mathfrak{C}\{L\} = -\frac{\log(1 + \lambda i)}{\log(1 - \lambda)},$$

while Proposition 2.6 yields the a priori bound

$$i\frac{1-\lambda}{1+i\lambda} \leq \mathfrak{C}\{L\} \leq i.$$

The function $\lambda \mapsto \mathfrak{C}$ decreases from $\mathfrak{C} = \mathfrak{C}_0 = i^-$ at $\lambda = 0^+$ to $\mathfrak{C} = 0^+$ at $\lambda = 1^-$. Hence, this function is discontinuous at $\lambda = 0$ since $\lambda = 0$ yields linear trees such that $\mathfrak{C} = 0$.

**4. A stochastic recursion.** To deal with the finiteness of $h$ in the general continuous case, we return to the recursion in discrete time. Recall that $L'$ denotes a random variable distributed like $L$, conditioned by $\{L \neq 0\}$, and that

$$\mathfrak{m}_1\{L\} := \mathbb{E}\{L \mid L \neq 0\}, \qquad \mathfrak{m}_2\{L\} := \mathbb{E}\{L^2\}/\mathbb{E}\{L\}.$$

Hence, for every $u$ in $(0, 1]$,

$$\mathfrak{m}_1\{L\} = \mathfrak{m}_1\{L'\} = \mathfrak{m}_1\{L_u\}, \qquad \mathfrak{m}_2\{L\} = \mathfrak{m}_2\{L'\} = \mathfrak{m}_2\{L_u\}.$$

4.1. *Results.* Our first aim in this section is to prove Proposition 4.1.

PROPOSITION 4.1. *Assume that $L$ is such that $\mathfrak{C}_0\{L\}$ is finite. Then $h(k)$ is finite if and only if $k > \mathfrak{m}_1\{L\}$.*

Corollary 4.2 follows directly and is anecdotal.

COROLLARY 4.2. *For any $L$ such that $\mathfrak{C}_0$ is finite, $h(k, b)$ is finite if and only if $k > \max\{\mathfrak{m}_1, b\}$.*

Our second aim is to characterize the sequence $(h(k))_k$ as solution of a stochastic recursion.

PROPOSITION 4.3. *For every $L$ such that $\mathfrak{C}_0\{L\}$ is finite and every $k > \mathfrak{m}_1\{L\}$,*

(2)  $$(k - \mathbb{E}\{L\})h(k) = k\mathbb{E}\{h(k+L)\}.$$

*This is equivalent to the fact that, for every $k > \mathfrak{m}_1\{L\}$,*

$$(k - \mathfrak{m}_1\{L\})h(k) = k\mathbb{E}\{h(k+L')\} = k\mathbb{E}\{h(k+L) \mid L \neq 0\}.$$

4.2. *Remark.* As was expected, the distributions $p$ and $p_u$ yield the same functional $h$. To check this a posteriori, assume that $h$ solves the recursion in Proposition 4.3 for $L$ of distribution $p$ and that $h^u$ solves the recursion in Proposition 4.3 for $L_u$. (Caution: $h^u$ corresponds to the continuous-time process with offspring distribution $p_u$ and $H^u$ to the discrete-time process $Z_u$.) Then $\mathbb{E}\{L_u\} = u\mathbb{E}\{L\} = um$ and Proposition 4.3 for $h^u$ becomes

$$(k - um)h^u(k) = k\mathbb{E}\{h^u(k + L_u)\}.$$

By definition of $L_u$,

$$\mathbb{E}\{h^u(k + L_u)\} = (1 - u)h^u(k) + u\mathbb{E}\{h^u(k + L)\}.$$

Rearranging and dividing both sides by $u$, one gets

$$(k - m)h^u(k) = k\mathbb{E}\{h^u(k + L)\}.$$

Since $h$ and $h^u$ satisfy the same relation, $h = h^u$ for every $u$ in $(0, 1]$. Finally, since the distributions of $L_u$ and $(L')_{u(1-p(0))}$ coincide, $L$ and $L'$ also yield the same function $h$.

4.3. *Proofs.*

PROOF OF PROPOSITION 4.1. Assume without loss of generality that $L \geq 1$ almost surely and let $m := \mathfrak{m}_1\{L\} = \mathbb{E}\{L\}$. Every $H_n^u(k)$ is finite, $H_0^u(k) = 1/k$, every sequence $(H_n^u(k))_n$ is nondecreasing and, for every $n \geq 1$,

(3) $$H_n^u(k) = (1 + um)\mathbb{E}_k\{H_{n-1}^u(Z_u(1))\}.$$

First, since $Z_u(1) = k$ with probability $(1 - u)^k$ on $\{Z_u(0) = k\}$,

$$H_n^u(k) \geq (1 + um)(1 - u)^k H_{n-1}^u(k).$$

If $k < m$, then $(1 + um)(1 - u)^k > 1$ for every $u$ sufficiently small. This implies that $H^u(k)$ is infinite for every $u$ sufficiently small, hence $h(k)$ is infinite.

Herefar, we assume that $k > m$. Since $H_{n-1}^u$ is nonincreasing, one has $H_{n-1}^u(Z_u(1)) \leq H_{n-1}^u(k + 1)$ on the event $\{Z_u(1) \geq k + 1\}$. The probability of $\{Z_u(1) \geq k + 1\}$ when $Z_u(0) = k$ is $1 - (1 - u)^k$, hence (3) implies that

$$H_n^u(k) \leq (1 + um)((1 - u)^k H_{n-1}^u(k) + (1 - (1 - u)^k)H_{n-1}^u(k + 1)).$$

Since $1 - (1 - u)^k \leq ku$ and $H_{n-1}^u \leq H_n^u$, we have

$$H_n^u(k) \leq (1 + um)((1 - u)^k H_n^u(k) + kuH_n^u(k + 1)).$$

Since $k > m$, for every positive and sufficiently small $u$, $(1 + um)(1 - u)^k < 1$. Hence, for $u$ sufficiently small,

$$H_n^u(k) \leq H_n^u(k + 1)\frac{ku(1 + um)}{1 - (1 + um)(1 - u)^k}.$$

Letting $n \to \infty$, the same inequality holds true for $H^u(k)$ and $H^u(k+1)$. When $u \to 0$, for every fixed $k > m$, the fraction on the right-hand side converges to the finite limit $k/(k-m)$, hence if $h(k+1)$ is finite, then $h(k)$ is finite too, and
$$h(k) \leq h(k+1)k/(k-m).$$
Finally, if $h(k)$ is finite for $k$ sufficiently large, and indeed $h(k)$ is finite for $k > \mathfrak{C}_0\{L\}$, then $h(k)$ is finite for every $k > m$.

The case $k = m$ is similar to the case $k < m$. For every $k$,
$$H_n^u(k) \geq (1+um)((1-u)^k H_n^u(k) + ku(1-u)^{k-1} H_n^u(k+1)).$$
One concludes by writing this as $H_n^u(k) \geq a_m^u(k) H_n^u(k+1)$ for an explicit $a_m^u(k)$ and by checking that $a_m^u(m) \to \infty$ when $u \to 0$. □

PROOF OF PROPOSITION 4.3. The sequence $(H^u(k))_k$ solves the recursion
$$H^u(k) = (1+um)\mathbb{E}_k\{H^u(Z_u(1))\}.$$
When $Z_u(0) = k$ and $u \to 0$, either (i) every individual survives, this happens with probability $1 - ku + o(u)$ and then $Z_u(1) = k$, or (ii) exactly one individual dies and is replaced by $1+L$ individuals, this happens with probability $ku + o(u)$ and then $Z_u(1) = k+L$. The probability that at least two individuals die is $o(u)$. Hence,
$$H^u(k) = (1+um)((1-ku)H^u(k) + ku\mathbb{E}\{H^u(k+L)\}) + o(u).$$
It follows that $h(k) = \lim_{u \to 0} H^u(k)$, if finite, solves the recursion stated in the proposition. The finiteness of $h(k)$ is crucial in this argument since one must cancel $h(k)$ on both sides of the equation, and this finiteness is guaranteed by Proposition 4.1. □

In relation to the last sentence of the proof of Proposition 4.3, note that if $\mathbb{E}\{L\} < k \leq \mathfrak{m}_1$, equation (2) in Proposition 4.3 is still valid in the sense that both sides are infinite.

**5. The first threshold.** Proposition 4.3 in Section 4.1 shows that $h(k)$ cannot be finite if $k \leq \mathfrak{m}_1$. Hereafter, we assume that $L \geq 1$ almost surely and that $k > \mathfrak{m}_1 = m := \mathbb{E}\{L\}$. In the special case $p = \delta_i$, $L = i$ almost surely and iterating Proposition 4.3 yields
$$(k-i)h(k) = (k+ji)h(k+(j+1)i) \geq \frac{k+ji}{k+(j+1)i}$$
for every integer $j \geq 0$, hence $h(k) \geq 1/(k-i)$. Since the reversed inequality was established above, the proof of Theorem 1 is complete.

Returning to the general case, we now prove part (b) of Theorem 2, which can be viewed as an extremality property of the deterministic branching mechanisms since it asserts that, for every other branching mechanism, for every $k > \mathfrak{m}_1\{L\}$,

$$h(k) > 1/(k - \mathfrak{m}_1\{L\}).$$

PROOF OF PART (b) OF THEOREM 2. One can assume without loss of generality that $L \geq 1$ almost surely. Let $m := \mathbb{E}\{L\}$. Since $h(k+L) \geq 1/(k+L)$, we have

$$(k-m)h(k) \geq k\mathbb{E}\{1/(k+L)\} \geq k/(k+m).$$

Thus, the a priori lower bound $h(k) \geq 1/k$ can be improved to the lower bound

$$h(k) \geq \frac{k}{(k-m)(k+m)}.$$

One can iterate this reasoning as follows. Assume that $h(k) \geq R(k)$, where $R(y)$ denotes a rational fraction, convex on $y > m$. Then

$$\mathbb{E}\{h(k+L)\} \geq \mathbb{E}\{R(k+L)\} \geq R(k+m),$$

and hence $h(k) \geq T(R)(k)$ for a new rational fraction $T(R)$, convex on $y > m$ and defined as

$$T(R)(y) := \frac{R(y+m)y}{y-m}.$$

Starting from $R(y) := 1/y$ and iterating this, one gets $h(k) \geq T^n(R)(k)$ with

$$T^n(R)(y) = \frac{1}{n+1}\left(\frac{1}{y+nm} + \frac{n}{y-m}\right).$$

Letting $n \to \infty$, this yields $h(k) \geq 1/(k-m)$. Finally, $h(k) > 1/(k-m)$ as soon as some strict convexity is involved, that is, for every distribution $p \neq \delta_m$. For instance, one can write $(k-m)h(k)$ as

$$k\mathbb{E}\{h(k+L)\} \geq k\mathbb{E}\{1/(k+L-m)\} > k/(k+\mathbb{E}\{L\}-m) = 1. \quad \square$$

**6. The second threshold.** Our aim in this section is to prove Theorem 5 stated in Section 2.3.

6.1. *Results.* It seems that the parameter $\mathfrak{m}_2$ introduced in Definition 1.5 enters the picture mainly because of the technical result in Lemma 6.2 below.

DEFINITION 6.1. For every positive $c$ and every $s$ in $[0,1]$, let

(4) $$\varphi_c(s) := c(1 - \mathbb{E}\{s^L\}) - \mathbb{E}\{L\}(1 - s^c).$$

LEMMA 6.2. *Let $c > 0$. Then $c \geq \mathfrak{m}_2$ if and only if $\varphi_c(s) \geq 0$ for every $s$ in the interval $[0, 1]$.*

LEMMA 6.3. *For any given $L$, $\mathfrak{C}_0\{L\} \geq \mathfrak{m}_2\{L\}$.*

LEMMA 6.4. *For every $L$, $\mathfrak{C}\{L\} \leq \mathfrak{m}_2\{L\}$. Since $\mathfrak{m}_2\{L'\} = \mathfrak{m}_2\{L\}$ and, for every $u$ in $(0, 1]$, $\mathfrak{m}_2\{L_u\} = \mathfrak{m}_2\{L\}$, this implies that $\mathfrak{C}\{L'\} \leq \mathfrak{m}_2\{L\}$ and $\mathfrak{C}\{L_u\} \leq \mathfrak{m}_2\{L\}$.*

This implies the upper bounds in Theorem 5.
Finally, we note that for every $L$ and every $u$ in $(0, 1]$,
$$\mathfrak{C}\{L'\} \leq \mathfrak{C}\{L\} \leq \mathfrak{C}\{L_u\} \leq \mathfrak{C}_0\{L\} = \mathfrak{m}_2\{M\},$$
where $M$ may be any of the random variables $L'$, $L$ or $L_u$.

### 6.2. Proofs.

PROOF OF LEMMA 6.2. First, $\varphi_c(s) = \psi_c(s^c)$ with
$$\psi_c(s) := c(1 - \mathbb{E}\{s^{L/c}\}) - \mathbb{E}\{L\}(1 - s).$$
Then $\psi_c(1) = \psi_c'(1) = 0$ and $\psi_c''(1) = \mathbb{E}\{L\}(1 - \mathfrak{m}_2/c)$. Hence, for any positive $c < \mathfrak{m}_2$, some values of $\varphi_c(s)$ are negative.

Conversely, assume that $c \geq \mathfrak{m}_2$. If $\psi_c$ is convex, then $\varphi_c \geq 0$ on $[0, 1]$. To show that $\psi_c$ is indeed convex, we first note that $\psi_c''(s^{1/c})$ has the same sign as $\mathbb{E}\{L(c - L)s^L\}$, which increases as $c$ increases. Thus, it is enough to show that this last expectation is nonnegative when $c = \mathfrak{m}_2$. In turn, this condition is equivalent to $\chi \geq 0$ on $[0, 1]$, where
$$\chi(s) := \mathbb{E}\{L^2\}\mathbb{E}\{Ls^L\} - \mathbb{E}\{L\}\mathbb{E}\{L^2 s^L\}.$$
One sees that $\chi(s) = \mathbb{E}\{K(s)\}/2$, where
$$K(s) := L_1 L_2 (L_1 - L_2)(s^{L_2} - s^{L_1})$$
and $L_1$ and $L_2$ are two independent copies of $L$. Since the function $L \mapsto s^L$ is nonincreasing, $K(s) \geq 0$ almost surely and hence $\chi(s) \geq 0$. This concludes the proof of Lemma 6.2. □

FIRST PROOF OF LEMMA 6.3. For every positive $c \geq \mathfrak{C}_0$, $L_u$ belongs to $\mathcal{G}_c$ for every $u$, in particular when $u \to 0$ [we recall that the family $(\mathcal{G}_c)_c$ is increasing]. Rearranging both sides of the inequality which characterizes the property that $L_u$ belongs to $\mathcal{G}_c$ in Definition 2.2 and taking logarithms, one sees that this inequality is equivalent to the fact that $\omega_L(s, u, c) \leq 0$ for every $s$ in $[0, 1]$ where

(5)  $\omega_L(s, u, c) := \log(1 + u\mathbb{E}\{L\}(1 - s^c)) + c\log(1 - u(1 - \mathbb{E}\{s^L\})).$

Since $\omega_L(s, u, c) = -\varphi_c(s)u + o(u)$ when $u \to 0$, one gets $\varphi_c(s) \geq 0$ for every $s$ in $[0, 1]$, that is, $c \geq \mathfrak{m}_2$; see Lemma 6.2. Finally, every $c \geq \mathfrak{C}_0$ is such that $c \geq \mathfrak{m}_2$, hence Lemma 6.3 holds. □

SECOND PROOF OF LEMMA 6.3. Theorem C in [7] states that for any $L$ in $\mathcal{G}_c$,

$$\sigma^2\{L\} \leq c\mathbb{E}\{L\}(\mathbb{E}\{L\} + 1).$$

We apply this to $L_u$. Since $\mathbb{E}\{L_u\} = u\mathbb{E}\{L\}$ and $\mathbb{E}\{L_u^2\} = u\mathbb{E}\{L^2\}$, this inequality is equivalent to $\mathfrak{C}\{L_u\} \geq \mathfrak{m}_2\{L\}(1 + o(1))$ when $u \to 0$, hence $\mathfrak{C}_0\{L\} \geq \mathfrak{m}_2\{L\}$. □

PROOF OF LEMMA 6.4. Using the notation from the proof of Lemma 6.3, we must show that $\omega_L(s, u, \mathfrak{m}_2) \leq 0$.

Since $\log(1 + v) \leq v$ for every $v$, $\omega_L(s, u, c) \leq \varsigma_L(s, u, c)$ for any $c$ with

$$\varsigma_L(s, u, c) := u\mathbb{E}\{L\}(1 - s^c) + c\log(1 - u(1 - \mathbb{E}\{s^L\})).$$

The first term in $\varsigma_L(s, \cdot, c)$ is linear with respect to $u$. The second term has the form $\beta \log(1 - \alpha u)$ with positive $\alpha$ and $\beta$, hence it is concave with respect to $u$. Thus, $\varsigma_L(s, \cdot, c)$ is concave. Since $\varsigma_L(s, 0, c) = 0$, we have $\varsigma_L(s, u, c) \leq 0$ for every $u$ in $[0, 1]$ if and only if the first derivative of $\varsigma_L(s, \cdot, c)$ at $0$ is nonpositive. Since this derivative is $-\varphi_c(s)$, Lemma 6.2 above shows that $\varsigma_L(\cdot, \cdot, c) \leq 0$ if and only if $c \geq \mathfrak{m}_2$. □

**7. Between the two thresholds.** Theorem 2 describes $h(k)$ when $k \leq \mathfrak{m}_1$ and when $k > \mathfrak{m}_2$. In this section, we study the regime $\mathfrak{m}_1 < k \leq \mathfrak{m}_2$. We start from the observation, drawn from Proposition 4.3, that for any $k > \mathfrak{m}_1$,

$$(k - \mathfrak{m}_1)h(k) \leq kh(k + 1).$$

For instance, if $\mathfrak{m}_2 - 1 < k \leq \mathfrak{m}_2$, then $k + 1 > \mathfrak{m}_2$ and

$$h(k + 1) \leq 1/(k + 1 - \mathfrak{m}_2).$$

Thus, for every $k > \mathfrak{m}_1$ such that $k + 1 > \mathfrak{m}_2$,

$$h(k) \leq k/((k - \mathfrak{m}_1)(k + 1 - \mathfrak{m}_2)).$$

One can iterate this trick to deal with values of $k > \mathfrak{m}_1$ such that $k + 2 > \mathfrak{m}_2$, and so on. This yields Proposition 7.1.

PROPOSITION 7.1. *Let $k > \mathfrak{m}_1$ and let $\ell$ denote any nonnegative integer such that $k + \ell > \mathfrak{m}_2$. Then*

$$h(k) \leq \frac{1}{k + \ell - \mathfrak{m}_2} \prod_{i=0}^{\ell-1} \left\{1 + \frac{\mathfrak{m}_1}{k + i - \mathfrak{m}_1}\right\}.$$

Every denominator in the upper bound for $h(k)$ is positive. A strange feature of this result is that the best upper bound for $h(k)$ does not always correspond to the smallest value of $\ell$ such that $k+\ell > \mathfrak{m}_2$. Simple computations show that the optimal choice is the smallest value of $\ell$ such that

$$k+\ell > \mathfrak{m}_1 \frac{\mathfrak{m}_2 - 1}{\mathfrak{m}_1 - 1} = \mathfrak{m}_2 + \frac{\mathfrak{m}_2 - \mathfrak{m}_1}{\mathfrak{m}_1 - 1} \geq \mathfrak{m}_2.$$

## 8. Branching processes and i.i.d. sums.

8.1. *Random products.* Our starting point in this section is the recursion of Proposition 4.3 which motivates the following definition:

DEFINITION 8.1. Let $X$ and $X_i$ denote i.i.d. nonnegative integrable random variables. Let $S_0 := 0$ and $S_n := X_1 + \cdots + X_n$ for every $n \geq 1$. For every positive $x$, let

$$R_n(x, X) := \prod_{i=0}^{n}\left\{1 + \frac{\mathbb{E}\{X\}}{x + S_i}\right\}, \qquad \varrho(x, X) := \limsup_{n \to \infty} \mathbb{E}\{R_n(x, X)/n\}.$$

We first deal with the finiteness of the functional $\varrho$, then we relate it to the harmonic setting. Lemma 8.2 lists some elementary facts, the proofs of which are omitted.

LEMMA 8.2. *Let $X$ denote a nonnegative, nonzero random variable.*

(i) *The function $x \mapsto \varrho(x, X)$ is nonincreasing.*
(ii) *For every nonnegative $n$, $\mathbb{E}\{R_n(x, X)\} \geq 1 + (n+1)\mathbb{E}\{X\}/x$.*
(iii) *Hence, $\varrho(x, X) \geq \mathbb{E}\{X\}/x$. In particular, $\varrho(x, X)$ is always positive.*
(iv) *For every positive $\lambda$, $\varrho(\lambda x, \lambda X) = \varrho(x, X)$.*

8.2. *Results.*

DEFINITION 8.3. For every nonnegative, not identically zero, integrable $X$, let

$$\mathfrak{m}_0\{X\} := \mathbb{P}\{X = 0\}\mathfrak{m}_1\{X\} = \mathfrak{m}_1\{X\} - \mathbb{E}\{X\}.$$

Note that $\mathfrak{m}_0\{X\}$ satisfies

$$(1 + \mathbb{E}\{X\}/\mathfrak{m}_0\{X\})\mathbb{P}\{X = 0\} = 1.$$

PROPOSITION 8.4. *For every nonnegative integrable $X$, $\varrho(x, X)$ is infinite for every $x \leq \mathfrak{m}_0\{X\}$.*

Proposition 8.5 implies a dichotomy result, stated as Corollary 8.6.

PROPOSITION 8.5.   *There exist finite positive functions $v(\cdot, X)$ and $w(\cdot, X)$ such that, for every positive $x > \mathfrak{m}_0\{X\}$ and every positive $y$,*

$$\varrho(x+y, X) \leq \varrho(x, X) \leq \varrho(x+y, X) v(x, X) w(x, X)^y.$$

COROLLARY 8.6.   *Either $\varrho(x, X)$ is infinite for every $x > \mathfrak{m}_0\{X\}$, or $\varrho(x, X)$ is positive and finite for every $x > \mathfrak{m}_0\{X\}$.*

Theorem 6 links some functions $\varrho$ to the harmonic moments of branching processes:

THEOREM 6.   *For every square integrable $L$ with nonnegative integer values and for every integer $k > \mathfrak{m}_1\{L\}$, the normalized limiting harmonic moment $h(k)$ which corresponds to the branching mechanism based on $L$ satisfies the relation*

$$h(k) = \frac{\varrho(k - \mathbb{E}\{L\}, L)}{\mathbb{E}\{L\}} = \frac{\varrho(k - \mathfrak{m}_1\{L\}, L')}{\mathfrak{m}_1\{L\}}.$$

Corollary 8.7 below is a direct consequence and complements Corollary 8.6.

COROLLARY 8.7.   *Let $X$ denote a nonnegative, integer, square integrable random variable. Then $\varrho(x, X)$ is finite if and only if $x > \mathfrak{m}_0\{X\}$. Furthermore, $x\varrho(x, X) \to \mathbb{E}\{X\}$ when $x \to \infty$. More precisely, for every integer $x > \mathfrak{m}_0\{X\}$,*

$$\mathbb{E}\{X\}/x \leq \varrho(x - \mathbb{E}\{X\}, X) \leq \mathbb{E}\{X\}/(x - \mathfrak{m}_2\{X\}).$$

From assertion (iv) in Lemma 8.2, the same conclusion holds if $X$ is a multiple of a nonnegative, integer, square integrable random variable. The hypothesis that $X$ is lattice and that $x$ belongs to this lattice might be unnecessary.

Note that one uses the square integrability of $X$ only to ensure the condition that $\limsup_{k\to\infty} kh(k) \leq 1$ [which is equivalent to $\lim_{k\to\infty} kh(k) = 1$].

For every positive $x$ and every nonnegative random variable $X$, the function $\varphi := \varrho(\cdot, X)$ solves the recursion

(6) $$x\varphi(x) = \mathbb{E}\{\varphi(x+X)\}(x + \mathbb{E}\{X\}).$$

Recall that $\varrho(x, X) \geq \mathbb{E}\{X\}/x$ in full generality. A consequence of (6) is that

$$x\varrho(x, X) \geq \mathbb{E}\{X\}\mathbb{E}\{x+X\}\mathbb{E}\{1/(x+X)\},$$

hence $\varrho(x, X) > \mathbb{E}\{X\}/x$ in full generality.

8.3. *Proofs.*

PROOF OF PROPOSITION 8.4. Let $r := \mathbb{P}\{X = 0\}$. Since $X_1$ is 0 with probability $r$ or distributed like $X'$ with probability $1 - r$, we have

(7)
$$\mathbb{E}\{R_{n+1}(x, X)\} = (1 + \mathbb{E}\{X\}/x)(r\mathbb{E}\{R_n(x, X)\} + (1 - r)\mathbb{E}\{R_n(x + X', X)\}).$$

The second part of Lemma 8.2 yields

$$\mathbb{E}\{R_n(x + X', X)\} \geq 1 + (n + 1)\mathbb{E}\{X\}\mathbb{E}\{1/(x + X')\}.$$

If $(1 + \mathbb{E}\{X\}/x)r \geq 1$, this implies that there exists a positive constant $c$ such that

$$\mathbb{E}\{R_{n+1}(x, X)\} \geq \mathbb{E}\{R_n(x, X)\} + c(n + 1).$$

Hence, $\mathbb{E}\{R_n(x, X)\} \geq cn^2/2$ and $\varrho(x, X)$ is infinite. Finally, the condition $(1 + \mathbb{E}\{X\}/x)r \geq 1$ can be rewritten as

$$x \leq \mathbb{E}\{X\}r/(1 - r) = \mathfrak{m}_0\{X\}.$$

This proves Proposition 8.4. □

PROOF OF PROPOSITION 8.5. Let $T$ denote the hitting time of the level $y$ by the nondecreasing process $(S_i)_{i \geq 0}$, that is,

$$T = \inf\{i \geq 0; S_i \geq y\}.$$

Since $S_i < y$ if $i \leq T - 1$ and $S_i \geq y$ if $i \geq T$,

$$\mathbb{E}\{R_n(x, X)|T\} \leq R_{T+n}(x, X) \leq R_{T-1}(x, X)R'_n(x + S_T, X)$$

and hence

$$\mathbb{E}\{R_n(x, X)|T\} \leq (1 + \mathbb{E}\{X\}/x)^T R'_n(x + y, X),$$

where $R'_n(\cdot, X)$ denotes a copy of $R_n(\cdot, X)$, independent of $T$. Hence,

$$\varrho(x, X) \leq \varrho(x + y, X)\mathbb{E}\{s^T\}, \qquad s := 1 + \mathbb{E}\{X\}/x.$$

To bound the last term, we use the fact that $\{T > n\} = \{S_n < y\}$ and hence an exponential inequality yields

$$\mathbb{P}\{T > n\} \leq u^{-y}\mathbb{E}\{u^X\}^n$$

for every positive $u \leq 1$. When $u \to 0$, $\mathbb{E}\{u^X\} \to \mathbb{P}\{X = 0\}$, hence the assumption on $x$ implies that there exists $u$ such that $s\mathbb{E}\{u^X\} < 1$. Hence,

$$\mathbb{E}\{s^T\} \leq \sum_{n \geq 1} s^n \mathbb{P}\{T > n - 1\} \leq u^{-b}\sum_{n \geq 0}(s\mathbb{E}\{u^X\})^n$$

and the last series converges. This proves the proposition. □

PROOF OF THEOREM 6. The choice $X = L'$ gives $\mathbb{E}\{X\} = \mathfrak{m}_1\{L\}$. Hence, Proposition 4.3 yields that, for every $k > \mathfrak{m}_1\{L\}$ and every $n \geq 0$,

$$h(k) = \mathbb{E}\{R_n(x, L')h(k + S_{n+1})\}, \qquad x := k - \mathfrak{m}_1\{L\}.$$

Since $R_n(x, L')$ and $h(k + S_{n+1})$ are nonincreasing functionals of the sequence $(X_i)_i$, a coupling inequality reads

$$h(k) \geq \mathbb{E}\{R_n(x, L')/n\}n\mathbb{E}\{h(k + S_{n+1})\}.$$

Now, by convexity,

$$\mathbb{E}\{h(k + S_{n+1})\} \geq \mathbb{E}\{1/(k + S_{n+1})\} \geq 1/(k + \mathfrak{m}_1\{L\}(n+1)).$$

When $n \to \infty$, we obtain $h(k) \geq \varrho(x, L')/\mathfrak{m}_1\{L\}$.

On the other hand, for every nonnegative integer $i$,

$$h(k) \leq \mathbb{E}\{R_n(x, L')\}h(k + i) + \mathbb{E}\{R_n(x, L'); S_n \leq i\}.$$

First, when $L'$ is square integrable, $h(k + i) \sim 1/i$ when $i \to \infty$. Second, we use the following special case of Cramér large deviations bounds:

LEMMA 8.8. *Let $X$ and $X_i$ denote i.i.d. nonnegative random variables. For every $y < \mathbb{E}\{X\}$, there exists $r < 1$ such that for every $n$,*

$$\mathbb{P}\{X_1 + \cdots + X_n \leq ny\} \leq r^n.$$

We apply this lemma to $X = L'$. This yields $\mathbb{P}\{S_n \leq ny\} \leq r^n$ with $r < 1$ for every $y < \mathfrak{m}_1\{L\}$. Third, since $L' \geq 1$, $R_n(x, L')$ is at most the product of $1 + \mathfrak{m}_1\{L\}/(x + j)$ from $j = 0$ to $n$, which is $O(n^{\mathfrak{m}_1\{L\}})$. Thus, for every $y < \mathfrak{m}_1\{L\}$, using the above for $i := \lfloor ny \rfloor$, one gets

$$h(k) \leq \varrho(x, L')(n/\lfloor ny \rfloor + o(1)) + O(n^{\mathfrak{m}_1\{L\}})r^n.$$

This suffices to prove Theorem 6. □

FIRST PROOF OF LEMMA 8.8. Let $X_i^t := \inf\{X_i, t\}$. For any $t$,

$$X_1^t + \cdots + X_n^t \leq X_1 + \cdots + X_n$$

and, if $t$ is sufficiently large, $y < \mathbb{E}\{X_1^t\}$. Since $X_1^t$ is almost surely bounded and the sequence $(X_n^t)_{n \geq 1}$ is i.i.d., standard large deviations bounds apply and hence there exists $r < 1$ with

$$\mathbb{P}\{X_1 + \cdots + X_n \leq ny\} \leq \mathbb{P}\{X_1^t + \cdots + X_n^t \leq ny\} \leq r^n. \qquad \square$$

A second proof of Lemma 8.8 uses the fact that for any nonnegative integrable random variable $X$, $\mathbb{E}\{e^{-tX}\} = 1 - t\mathbb{E}\{X\} + o(t)$ when $t$ is positive and $t \to 0$.

SECOND PROOF OF LEMMA 8.8. Assume first that $X$ is integrable. By Chebyshev's inequality, the result of the lemma holds with

$$r := e^{ty}\mathbb{E}\{e^{-tX}\}$$

for every positive $t$. Since $X$ is integrable and positive, its Laplace transform $\mathbb{E}\{e^{-tX}\}$ is differentiable at $t = 0^+$, that is, $\mathbb{E}\{e^{-tX}\} = 1 - t\mathbb{E}\{X\} + o(t)$ when $t$ is positive and $t \to 0$. Since $e^{ty} = 1 + ty + o(t)$ and $y < \mathbb{E}\{X\}$, one obtains a value $r < 1$ for any $t$ positive and sufficiently small.

If $X$ is not integrable, truncating $X$ and every $X_n$ as in our first proof yields the result. $\square$

INSTITUT FOURIER
100 RUE DES MATHS
BP 74
38402 SAINT MARTIN D'HÈRES
FRANCE
E-MAIL: Didier.Piau@ujf-grenoble.fr
URL: www-fourier.ujf-grenoble.fr